\documentclass{article}

\usepackage{amsthm}
\theoremstyle{plain}
\newtheorem{thm}{Theorem}
\newtheorem*{thm*}{Theorem}
\newtheorem{lem}{Lemma}
\newtheorem{pro}{Proposition}
\newtheorem{cor}{Corollary}

\theoremstyle{definition}
\newtheorem{dfn}{Definition}
\newtheorem*{pf}{Proof}
\newtheorem{ex}{Example}

\usepackage{amsmath, amsfonts, amsthm, amssymb, amscd}
\usepackage{type1cm}
\usepackage{amscd}
\usepackage{mathrsfs}
\usepackage{graphicx}
\usepackage[mathscr]{euscript}
\usepackage{mathrsfs}
\usepackage{caption}
\usepackage{setspace}

\def\MARU#1{{\rm\ooalign{\hfil\lower.168ex\hbox{#1}\hfil \crcr\mathhexbox20D}}}

\begin{document}

\title{Dynamical systems for eigenvalue problems of axisymmetric matrices with positive eigenvalues}

\author{Shintaro Yoshizawa \thanks{R-Frontier Div. Frontier Research Center, Toyota Motor Corporation, 543, Kirigahora, Nishihirose, Toyota, Aichi, 470-0309 Japan. \quad  ${\rm shintaro\_yoshizawa@mail.toyota.co.jp}$  }}

\maketitle



\begin{abstract}
We consider the eigenvalues and eigenvectors of an axisymmetric matrix$A$ with some special structures.
We propose S-Oja-Brockett equation 
$\frac{dX}{dt}=AXB-XBX^TSAX,$
where $X(t) \in {\mathbb R}^{n \times m}$ with $m \leq n$, 
$S$ is a positive definite symmetric solution of the Sylvester equation 
 $A^TS = SA$ 
and 
$B$ is a real positive definite diagonal matrix whose diagonal elements are distinct each other, 
and show the S-Oja-Brockett equation has the global convergence to eigenvalues and its eigenvectors of $A$.
\end{abstract} 

\bigskip

\section{Introduction}
In least squares optimization, 
Brockett \cite{Bro1},\cite{Bro2},\cite{Bro3} showed that the tasks of diagonalizing a matrix, linear programming, and
sorting, could all be solved by dynamical systems given by
$$\frac{dX}{dt}=AXB-XBX^TAX,$$
where $X(t)$ belongs to the real special orthogonal group, that is, $X^TX=I$ and $\det(X)=1,$ and $A$, $B$ are real symmetric matrices. The symbol $T$ denotes the transpose of the matrix.
His results had their origins in earlier work
dating back to that of Fischer \cite{Fis}, Courant \cite{Cou} and von Neumann \cite{VonN}, 
Also there were parallel efforts in numerical analysis by Chu \cite{Chu}.

On the other hand, in the field of neural networks, 
Amari \cite{Ama} pointed out that the maximum eigenvalue problem can be obtained by the Hebbian learning rule, 
and then Oja \cite{Oja1},\cite{Oja2},\cite{Oja3} in a more generalized form, posited that the principal subspace can be obtained in dynamical systems. For a real positive definite symmetric matrix $A$, by linear transformation 
$Y = A^{1/2}X$ of the Oja equation given by
$$\frac{dX}{dt}=(I-XX^T)AX,$$
where $I$ is the identity matrix and $X(t)$ belongs to ${\mathbb R}^{n \times m}$ which denotes an n-by-m real matrix space. 
Wyatt-Elfadel \cite{WE} showed that the following equation is a gradient flow: 
$$\frac{dY}{dt}=AY-YY^TY, $$ 
However, Wyatt-Elfadel \cite{WE} did not realize that the gradient is defined for the Riemannian metric.
Xu \cite{Xu} proposed the Oja-Brockett equation given by
$$
\frac{dX}{dt}=AXB-XBX^TAX,
$$
where $X(t) \in {\mathbb R}^{n \times m}$, $A$ is a positive diagonal matrix and 
$B$ is a real positive definite symmetric matrix. 
Thus, Oja-Brockett equation is a generalization of the dynamical system proposed by Oja \cite{Oja3} and Brockett \cite{Bro1}, and investigated its properties, but it was believed at the time that the Oja-Brockett equation was not gradient flow. 
Local stability analysis (Xu \cite{Xu}) was developed near the equilibrium point, but it does not actually prove global convergence to the equilibrium point.

In contrast, Yoshizawa-Helmke-Starkov \cite{YHS} proved that the Oja-Brockett equation is a gradient flow with suitable Riemannian metric 
for a real positive definite symmetric matrix $A$, and that if the initial matrix $X(0)$ has full rank, the solution of the Oja-Brockett equation converges to eigenvectors of $A$ globally.This is
exactly what we are doing and our main result concerning global convergence is then
deduced using a result by \L ojasiewicz \cite{Lo} on real analytic gradient flows.
Of course, this implies that the Oja equation is a gradient flow, too. 

If the coefficient matrix $A$ of the Oja-Brockett equation is a non-symmetric matrix with positive eigenvalues, it unfortunately does not converge to the eigenvectors of $A$. 
Therefore, the purpose of this paper is to propose a new dynamical system(we call this S-Oja-Brockett equation) given by
\begin{equation}\label{sob}
\frac{dX}{dt}=AXB-XBX^TSAX, 
\end{equation}
where 
$X(t) \in {\mathbb R}^{n \times m}$, 
$B$ is a real positive definite symmetric matrix and  
$S$ is a real positive definite symmetric solution of the Sylvester equation given by 
\begin{equation}\label{syl}
A^TS = SA.
\end{equation}

We consider the eigenvalues and eigenvectors of an asymmetric matrix $A$ with some special structures, 
and show the equation (\ref{sob}) has the global convergence to eigenvectors and eigenvalue sorting properties of $A$ that the Oja-Brockett  equation has.

That is, the equation (\ref{sob}) is valid for a positive definite symmetric matrix $B$ and a coefficient matrix $A$ of the Sylvester equation (\ref{syl}) for which there exists a positive definite symmetric matrix solution $S$. 
As a special case, if $A$ is a positive definite symmetric matrix, then there is obviously the identity matrix as $S$.
The equation (\ref{sob}) is characterized by the ability to simultaneously obtain eigenvalues and eigenvectors globally, and also to control the extraction of eigenvalues and eigenvectors of $A$ by the magnitude of the eigenvalues of $B$.

\bigskip
This paper is organized as follows, In Section 2,
we describes the relationship between the solution $S$ and the coefficient matrix $A$ of the Sylvester equation (\ref{syl}).
In Section 3, we introduce the equation (\ref{sob}) and discuss 
an example of the equation (\ref{sob}) for the positive definite diagonal solution $S$ of the Sylvester equation (\ref{syl}).
In Section 4, we investigate another example of 
the equation (\ref{sob}) for 
the blocked positive definite symmetric solution $S$ of 
the Sylvester equation (\ref{syl}).

\section{Representation of the Sylvester equation}
In this section, we discuss the structure of the relationship between solutions and coefficients of the Sylvester equation.
The results of Taussky-Zassenhaus \cite{TZ} considered matrices that transform a square matrix into its transpose matrix.
We consider the existence and structure of a matrix that transforms a square matrix into a symmetric matrix in a different way from that of \cite{TZ}.
We prepare a Lemma to prove the following theorem.
\begin{thm}\label{theorem1}
A real regular matrix $A$ with distinct positive eigenvalues 
becomes a symmetric matrix by a similarity transformation by some regular matrices.
\end{thm}
\quad
\begin{lem}[\bf{Hodge \cite{Hodge}}]\label{lem1}
Let $A$ be a square matrix.
The Sylvester equation $A^{T}S=S^{T}A$ has solution $S$ and if $A$ is regular, then the general solution is given by $S = P^{T}ZQ^{-1},$ 
where $Z$ is constrained only by the symmetry requirement, that is, $Z^{T}=Z$ and $PAQ=I$.
\end{lem}
\begin{pf}\rm{ By direct computation. 
{\hfill $\Box$}
}\end{pf}  
\quad
\begin{pf}[ {\bf Theorem\ref{theorem1}} ] \rm{ 
Since $A$ has $n$ orthogonal eigenvectors, there is an orthogonal matrix $U$ and a positive diagonal matrix $D$ 
such that $U^{T}AU = D$ and 
$D^{-\frac{1}{2}}U^{T}AUD^{-\frac{1}{2}}=I$. Thus, we define $P=D^{-\frac{1}{2}}U^{T}$ and $Q=UD^{-\frac{1}{2}}$.
If we choose $Z$ a positive diagonal matrix, due to Lemma\ref{lem1}, 
we see that
$S = P^TZQ^{-1}$ becomes symmetric and has positive eigenvalues. Therefore, $S$ has the square root $S^{\frac{1}{2}}$.
Then $S^{\frac{1}{2}}AS^{-\frac{1}{2}}$ becomes a symmetric matrix since
$(S^{\frac{1}{2}}AS^{-\frac{1}{2}})^T =S^{-\frac{1}{2}} (S^{-\frac{1}{2}}SA)^T = S^{-\frac{1}{2}} (S^{-\frac{1}{2}}A^TS)^T = S^{\frac{1}{2}}AS^{-\frac{1}{2}}.$
Thus, we obtain the result.
{\hfill $\Box$}
}\end{pf} 
\quad
\begin{ex}
Consider the following matrix given in Tanabe-Sagae \cite{TS2},\cite{TS}:
$$
A = D + ab^T  \in \mathbb{R}^{n \times n}, \quad 
D = {\rm{diag}}(d_1,d_2,...,d_n) \in \mathbb{R}^{n \times n},  $$
$$
a = (a_1, a_2,...,a_n)^T , b = ( b_1, b_2,..., b_n )^T \in \mathbb{R}^{n}.
$$
We assume $0 < d_1 < d_2 <\cdots < d_n $ and $0 < a_ib_i,  (i = 1,2,...,n)$.
The characteristic polynomial of $A$, denoted by $\phi_{A}(\lambda) =\det(\lambda I - A)$, 
is given by 
$$\phi_{A}(\lambda) = \bigg(1 - \displaystyle\sum_{i=1}^{n}\frac{a_i b_i}{\lambda - d_i} \bigg)
\displaystyle\prod_{i = 1}^{n} ( \lambda - d_i ).
$$
Thus, we see $0 < d_1 < \lambda_1 < d_2 < \lambda_2 < d_3 < \cdots < d_n < \lambda_n.$
By $S =  {\rm{diag}}({\frac{b_1}{a_1}}, {\frac{b_2}{a_2}},\cdots, {\frac{b_n}{a_n}}
),$
$S^{\frac{1}{2}}AS^{-\frac{1}{2}} $ is a symmetric matrix.

\bigskip
We have a following commutative diagram:
\quad
$$ 
\begin{CD}
\xi = (d,(a,b)) \in \mathbb{R}^{n}_{+} \times \bigg( \mathbb{R}^{n} \setminus\{0\} \times \mathbb{R}^{n} \setminus\{0\} \bigg) 
@>{\varphi}_{\rm Mat}>> \mathbb{R}^{n \times n} \setminus\{0\} \ni \varphi_{\rm Mat}(\xi)   \\
@VV{\varphi}_{\rm Vec} V
@VV{\tilde{\pi}=\pi \circ \rm{Vec}}V \\
\varphi_{\rm Vec}(\xi) \in \mathbb{R}^{n^2} \setminus\{0\}
@>\pi>> \qquad \qquad \mathbb{P}^{n^2-1} 
\ni \eta \\
\end{CD}
$$
\quad 
where $\mathbb{P}^{n^2-1}$ is the real projective space of dimension $n^2-1$.  Let $\{e_i\}_{i=1}^{n}$ be the standard unit vectors and $E_i$ be the matrix with only the i-th diagonal 1. The maps in the diagram are given by 
$\varphi_{\rm Mat}(\xi) = D+ab^T, \quad \varphi_{\rm Vec}(\xi) = (\sum_{i=1}^{n} e_i \otimes E_i )d +b \otimes a$, 
\quad ${\rm Vec}(\varphi_{\rm Mat}(\xi)) = \varphi_{\rm Vec}(\xi)$ and 
$$
\tilde{\pi}(\varphi_{\rm Mat}(\xi)) 
= \pi(\varphi_{\rm Vec}(\xi)) =
\eta = \frac{\bigg\{\bigg( \displaystyle\sum_{i=1}^{n} e_i \otimes E_i \bigg)d + b \otimes a \bigg\}
\bigg\{\bigg( \displaystyle\sum_{i=1}^{n} e_i \otimes E_i \bigg)d + b \otimes a \bigg\} ^T}
{\bigg\{\bigg( \displaystyle\sum_{i=1}^{n} e_i \otimes E_i \bigg)d + b \otimes a \bigg\}^T
\bigg\{\bigg( \displaystyle\sum_{i=1}^{n} e_i \otimes E_i \bigg)d + b \otimes a \bigg\}}.
$$
\end{ex}

\section{Axisymmetric structure matrix}
The Oja-Brockett equation has the property that large eigenvalues and eigenvectors of a positive definite symmetric matrix can be simultaneously obtained (Yoshizawa-Helmke-Starkov \cite{YHS}), but surprisingly, by changing the coefficients of its equation, it is also find to be valid for a certain class of asymmetric matrices with real positive eigenvalues. 

In this section, as an extension of the concept of an band matrix, we consider axisymmetric structure matrices as coefficient matrices for positive definite diagonal matrix solutions of the Sylvester equation, and propose the eigenvalues of the axisymmetric structure matrix with real positive eigenvalues are obtained by the S-Oja-Brockett equation (\ref{sob}).

\begin{dfn}[{\bf Axisymmetric Structure Matrix}]\label{asm}  The axisymmetric structure matrix is definde by 
\begin{align*}
\mathscr{A}_n = \{ A \in \mathbb{R}^{n \times n} | A & = {\rm diag}(d_1,d_2,...,d_n) + (a_{ij}) + (b_{ji}), (j < i), 
\quad 0 \leq a_{ij} \cdot b_{ji},  (i \ne j),  \\
& {\rm if} \quad  a_{ij} \cdot b_{ji} = 0 \Rightarrow  a_{ij} = b_{ji} =0, \\ 
& {\rm s.t.} \quad [ \quad \exists S ={\rm diag}( s_1, s_2,..., s_n) >0 \smallskip  \quad {\rm s.t.} \quad  A^TS =SA  \quad ] \quad  \}.
\end{align*}
\end{dfn}

\bigskip
The axisymmetric structure matrix can be expressed in the form where the diagonal matrix is perturbed in rank as follows:

Let $A = D + \bigtriangleup_a + \bigtriangleup_b$, where $D = {\rm diag}(d_1,d_2,...,d_n), 
\bigtriangleup_a = (a_{ij}), \bigtriangleup_b = (b_{ji}), j<i.$ \\
If ${\mathbf{a}}_j = (0,...,0,1,a_{j+1 j},...,a_{n j})^T \in \mathbb{R}^{n}$ and 
${\mathbf{b}}_j = (0,...,0,1,b_{j j+1},...,b_{j n})^T \in \mathbb{R}^{n}$, then 
$\mathscr{A}_n$ can be expressed as a rank perturbation of the diagonal matrix as follows: 
\begin{equation}\label{rep}
A =  D + \displaystyle\sum_{j=1}^{n}({\mathbf{a}}_j {\rm e}_j^T + {\rm e}_j{\mathbf{b}}^T_j),
\end{equation}
where $D = {\rm diag}(d_1- 2, d_2 - 2, ..., d_n - 2)$ and the unit vector with kth element only 1, $e_k = (0,...,0,1,0,...,0)^T$.   \\
\begin{figure}[h]
  \centering
  \includegraphics[width=0.3\textwidth]{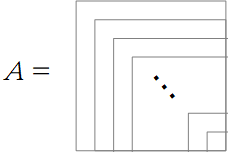}
  \caption{Subspace embedding of the axisymmetric structure matrix }\label{figure_1}
\end{figure}
Almost obvious, but leave it as a proposition.
\begin{pro}\label{emb}
The subspace generated by 
${\mathbf{a}}_j {\rm e}_j^T + {\rm e}_j{\mathbf{b}}^T_j$ is invariant under the similarity transformation of the regular diagonal matrix.
\end{pro}

${\mathbf{a}}_j {\rm e}_j^T + {\rm e}_j{\mathbf{b}}^T_j$ forms a band as shown in Figure \ref{figure_1}. If we restrict the solution $S$ of $A^TS=SA$ to be a symmetric matrix, 
then we do not lose generality by assuming $S$ to be a diagonal matrix.
In particular, we are interested in the positive definite solution $S$ 
to make $A\in \mathscr{A}_n$ a symmetric matrix by a similarity transformation.

In order to characterize the coefficient-solution pairs $(A,S)$ of equation $A^TS=SA$, 
we define the Lagrangian subspace of $\mathbb{R}^{2n}$ that we mean as follows.

\begin{dfn}\label{Lag}
We say ${\cal L} \subseteqq {\mathbb{R}^{2n}}$
is a Lagrangian subspace if $\cal L$ has dimension $n$ and
$$\langle J_n x, y \rangle_{\mathbb{R}^{2n}} = 
\langle x^1, y^2 \rangle_{{\mathbb R}^n} - \langle x^2, y^1 \rangle_{{\mathbb R}^n} = 0,$$
for all 
$x=
\left(
\begin{matrix}
x^1   \\
x^2
\end{matrix}
\right),
y=
\left(
\begin{matrix}
y^1   \\
y^2
\end{matrix}
\right)
\in {\cal L}$, 
$(x^i,y^i \in {\mathbb R}^n).$
Here, $\langle \cdot, \cdot \rangle_{\mathbb{R}^{2n}}, 
(\langle \cdot, \cdot \rangle_{\mathbb{R}^{n}})$ denotes Euclidean inner product on ${\mathbb{R}^{2n}}, ({\mathbb{R}^{n}} )$
, and
$$J_n
 =
\left(
\begin{matrix}
0 & -I_n  \\
I_n & 0 
\end{matrix}
\right),$$
with $I_n$ the $n$ by $n$ identity matrix.
\end{dfn}

More generally, any Lagrangian subspace of ${\mathbb{R}^{2n}}$
can be spanned by a choice of n linearly
independent vectors in ${\mathbb{R}^{2n}}$. 
We may regard it these $n$ vectors as the columns of a $2n$ by $n$ matrix$\tilde X$, 
which we shall call a frame for $\cal L$. 
Moreover, we will write 
$${\tilde X} =
\left(
\begin{matrix}
X \\
Y
\end{matrix}
\right),$$
where $X$ and $Y$ are $n$ by $n$ matrices.

\begin{pro}
Let $A \in \mathscr{A}_n$ and $S={\rm diag}(s_1,s_2,...,s_n)$. 
Then ${\tilde A}=\left(\genfrac{}{}{0pt}{}{A}{S}\right)$
is a frame for a Lagrangian subspace of ${\mathbb{R}^{2n}}$.

\end{pro}

\begin{pf}\rm{
Let $a_i$ and ${\hat s}_i$ be the i-th column vectors of matrices $A$ and $S$, respectively.
Since 
$0 = \left(A^TS-SA\right)_{ij} = 
\langle a_i, {\hat s}_j \rangle_{\mathbb{R}^{n}}
-\langle a_j, {\hat s}_i \rangle_{\mathbb{R}^{n}}
=\langle J
\left(\genfrac{}{}{0pt}{}{a_i}{{\hat s}_i}\right), \left(\genfrac{}{}{0pt}{}{a_j}{{\hat s}_j}\right)
\rangle_{\mathbb{R}^{2n}} ,$
we obtain the result.
Here, $()_{ij}$ denotes the ij component of the matrix.
{\hfill $\Box$}
}\end{pf}  

We define a {\it coordinate graph} for a square matrix to study the relationship 
between the solution and coefficients of the Sylvester equation, 
that is, we clarify the properties of the frame of the Lagrangian subspace.

\begin{dfn}[{\bf Coordinate Graph}]\label{gra}
For a matrix $A = (a_{ij})$, if $a_{ij} \neq 0$, then a {\it graph} for the matrix $A$ 
is defined by connecting $i$ and $j$ by edges, with row subscript $i$ and column subscript $j$ as nodes.
\end{dfn}
\begin{figure}[h]
\begin{center}
\includegraphics[width=1.0\textwidth]{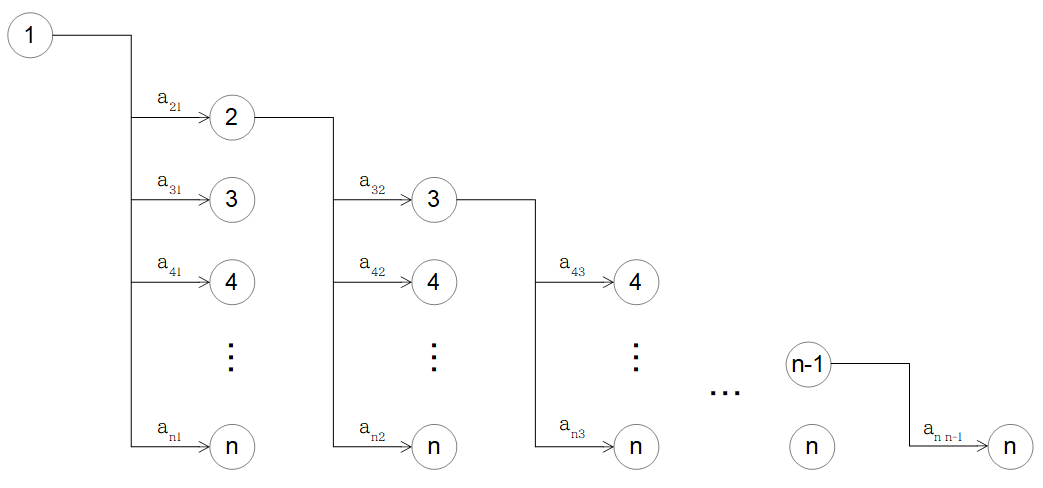}
\end{center}
\caption{Coordinate graph of the lower triangular part of matrix $A$ }\label{figure_2}
\end{figure}
Figure \ref{figure_2} shows the coordinate graph when all elements of the lower triangular part of matrix $A$ are non-zero. As can be seen from the graph, restricting the solution of the Sylvester equation to a diagonal matrix generally results in an over-determined system.

\begin{thm}[{\bf Connection equation}]\label{con} 
For the coordinate graph of the lower triangular part of matrix $A\in \mathscr{A}_n$, consider a family of matrices $A$ that can be defined as 
\begin{equation}\label{connection}
\frac{b_{j_{p}i}}{a_{ij_{p}}}\cdot s_{j_{p}} = \frac{b_{j_{q}i}}{a_{ij_{q}}}\cdot s_{j_{q}}, \quad 
( \quad j_{p}, j_{q} \in \{j_1,...,j_{r_i}\} \quad )
\end{equation}  
with respect to the equation 
$s_i=\frac{b_{ji}}{a_{ij}}\cdot s_j, \quad (\quad j \in \{j_1,...,j_{r_i}\} \subset \{1,2.3,...,i-1\} \quad )$ determined from the i-row component of the lower triangle of matrix $A$. If a representative node $s_*$ is arbitrarily chosen for each connected component of the graph and $s_*=1$, 
the positive definite diagonal matrix solution of the Sylvester equation can be constructed only by the lower and upper triangular parts, independent of the diagonal components of matrix $A$, and can be uniquely determined regardless of the coordinate system of the lower and upper triangular parts.
Here, $s_*$ is identified with $*$.
\end{thm}

\quad
\begin{pf}[ {\bf Theorem\ref{con}} ] \rm{ 
In Figure \ref{figure_2}, $\MARU{1}$ represents the seed for constructing the solution of the Sylvester equation, i.e., $s_1=1$. 
An edge is defined when a component of the lower triangular part of matrix $A$ exists.
Sequentially creating edges from the nodes in the left column to the right generates a tree structure.
If the tree structure, including the case of a single node, is not connected, one seed is set for each connected component.
At this time, if there are multiple elements in the same row of $A$, i.e., the same node,
In order to determine the solution of the Sylvester equation independently of the coordinate system of matrix $A$, 
by equating node number $j$ with $s_j$, we can equate the representation of $s_j$ for each connected component. 
In this way, thus, we obtain the conclusion.
{\hfill $\Box$}
}\end{pf} 
\quad

The equation (\ref{connection}) consisting of the constraint relation of the coefficients is called the {\it connection equation}.
We consider two concrete examples of the above theorem.

\quad
\begin{ex}[{\bf Non simply connected tree}]\label{tree1} 
Let $A_1 \in \mathscr{A}_n$ such that 
$$
A_1=
\begin{pmatrix}
d_1     &  0       &   0      &   0       &    0      \\
0         & d_2    & b_{23} &  0        & b_{25}  \\
0         & a_{32} & d_3    & b_{34}  & b_{35} \\
0         & 0        & a_{43} &  d_4    &  0        \\
0         & a_{52} & a_{53} &  0        & d_5  
\end{pmatrix}. 
$$
Then the graph is given in Figure \ref{figure_3} .
\begin{figure}[h]
\centering
\includegraphics[scale=0.5]{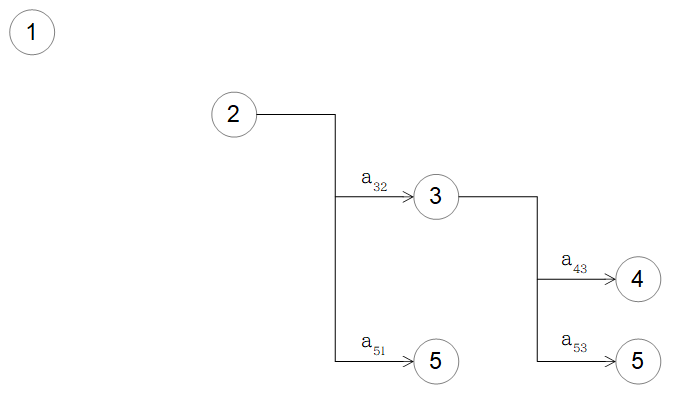}
\caption{Non simply connected tree}\label{figure_3}
\end{figure}

Since there are two connected components of the graph, $s_3, s_4$, and $s_5$ can be determined 
by $a_{ij} \cdot s_i = b_{ji} \cdot s_j$  if $s_1$ and $s_2$ are set to $s_1=1$ and $s_2=1$ as seeds, respectively. 
However, since $a_{52}$ and $a_{53}$ are in the same row of $A_1$, i.e., 
there are two ways to express ${s_5}$, we can uniquely determine the solution of the Sylvester equation if matrix $A_1$ satisfies the relation of the connection equation 
$\frac{b_{25}}{a_{52}} \cdot s_2 = \frac{b_{35}}{a_{53}} \cdot s_3$.
\end{ex}
\quad

\begin{ex}[{\bf Simply connected tree}]\label{tree2} 
Let $A_2 \in \mathscr{A}_n$ such that 
$$
A_2=
\begin{pmatrix}
d_1     &   0      & b_{13} & b_{14}  & b_{15} \\
0         & d_2    & b_{23} & 0         & b_{25}  \\
a_{31}  & a_{32} & d_3    & b_{34}  &  0        \\
a_{41}  &  0       & a_{43} &  d_4    & b_{45} \\
a_{51}  & a_{52} & 0        &  a_{54} & d_5  
\end{pmatrix}. 
$$ 
Then the graph is given in Figure \ref{figure_4}.
\begin{figure}[h]
\centering
\includegraphics[scale=0.5]{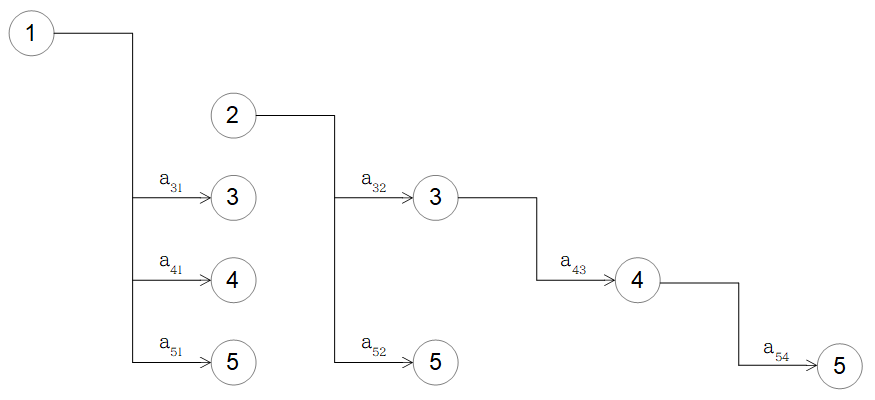}
\caption{Simply connected tree}\label{figure_4}
\end{figure}

Since $a_{31}$ and $a_{32}$, $a_{41}$ and $a_{43}$, as well as $a_{51}, a_{52}$, and $a_{54}$ 
are in the same row of $A_2$, respectively, the connection equations are defined as 
$\frac{b_{13}}{a_{31}}\cdot s_1=\frac{b_{23}}{a_{32}} \cdot s_2, 
\frac{b_{14}}{a_{41}}\cdot s_1=\frac{b_{34}}{a_{43}}\cdot s_3$ and 
$\frac{b_{15}}{a_{51}}\cdot s_1=\frac{b_{25}}{a_{52}}\cdot s_2=\frac{b_{45}}{a_{54}}\cdot s_4$. 
In this case, since the graph of $A_2$ is singly connected, the solution 
$S={\rm diag}(s_1,\dots,s_5)$  
of the Sylvester equation is uniquely determined by taking $s_1$ as the seed and setting $s_1=1$.
\end{ex}
\bigskip

We prepare the following complement to determine the sufficient conditions for 
$A \in \mathscr{A}_n$ to have positive eigenvalues.
\begin{lem}
The eigenvalues $\lambda$ of $A \in \mathscr{A}_n$ are contained in the following set:
\begin{equation}\label{eig}
\displaystyle\bigcup_{i=1,...,n}\{ \quad \lambda \in \mathbb{R} \quad \vert \quad
| \lambda - d_i | \quad \leq \displaystyle\sum_{j=1,...,n (j \ne i)} 
\sqrt{a_{ij} \cdot b_{ji}}  \quad \}.
\end{equation}
\end{lem}
\begin{pf}\rm{ By direct computations.
{\hfill $\Box$}
}\end{pf}  

\begin{cor}\label{cor1} 
Suppose the diagonal components of matrix $A\in \mathscr{A}_n$ satisfy 
$$\sum_{j=1,...,n  (j \ne i)} \sqrt{a_{ij}\cdot b_{ji}} \quad  < \quad  d_{i},  (i = 1,...,n),$$
then the eigenvalues of $A$ are positive.
\end{cor}
\begin{thm}[{\bf S-Oja-Brockett equation}]\label{SOB} 
For a matrix $A \in \mathscr{A}_n$ satisfying Corollary\ref{cor1}, 
we suppose that the initial values $X(0) \in \mathbb{R}^{n \times m}$ with $m \leq n$ are full rank and
the eigenvalues of matrix $A$ and the eigenvalues of positive definite symmetric matrix $B$ are distinct each other.
Then, we have 
\begin{flushleft}
{\rm (i)} The equation defined by 
\begin{equation}
\frac{dX}{dt} = AXB - XBX^TSAX
\end{equation}
is a real analytic gradient flow, and 
the matrix $X(t)^TSAX(t)$ converges to the $m$ largest eigenvalue matrix ${\rm diag}(\lambda_1,...,\lambda_m)$ of $A$ as $t \to \infty$, and at the same time 
the column vector $x_j$ of $X$ converges to the eigenvector of eigenvalue $\lambda_j$.  
Where the matrices $A, B$ and $S$ are time constants, and the matrix $S$ is a positive definite solution of the Sylvester equation.
\end{flushleft}
\begin{flushleft}
{\rm (ii)} If the matrix $B$ is represented as a diagonal matrix $B={\rm diag}(b_1,...,b_m)$ 
with $b_m < \cdots < b_1$, then
$X(t)^TSAX(t)$ converges to the matrix ${\rm diag}(\lambda_1,...,\lambda_m)$ 
with $\lambda_m < \cdots < \lambda_1$ and 
the i-th column vector of $\sqrt{S}X(\infty)$ is the eigenvector corresponding to the i-th eigenvalue $\lambda_i$ of $A$.
\end{flushleft}
\begin{flushleft}
{\rm (iii)} If matrix $B$ is simply a positive definite diagonal matrix, 
then in general, $X^TSAX$ 
does not converge to a diagonal matrix at time infinity, 
but $X$ converges to a matrix consisting of column vectors that generate the principal $m$ subspaces of $A$.
\end{flushleft}
\end{thm}
\begin{pf}\rm{
First, we show the S-Oja-Brockett equation is a real analytic gradient flow.
As we know in Yoshizawa-Helmke-Starkov \cite{YHS}, the flow given as 
$$\frac{dY}{dt}={\tilde A}YB-YBY^T{\tilde A}Y    ,  \quad Y \in R^{n \times m},$$ 
where $0 < \tilde{A}={\tilde{A}}^T $, $0 < B=B^T $, is a real analytic negative gradient flow with the potential function defined by 
$$f(Y)=\frac{1}{4}{\rm tr}[({\tilde A}YBY^T)^2]-\frac{1}{2}{\rm tr}({\tilde A}^2YB^2Y^T)$$ 
and the Riemannian metric defined by $$\langle \Omega_1, \Omega_2 \rangle ={\rm tr}(\tilde{A} \Omega_1 B \Omega_2^T), 
\quad \Omega_1,\Omega_2 \in {\mathbb R}^{n \times  m}.$$ 
That is, $\frac{dY}{dt}=-{\rm grad}f(Y)$. 
Since any matrix $A$ satisfying Corollary\ref{cor1} has a positive definite symmetric matrix solution $S$ of the Sylvester equation $A^TS=SA$, we can consider the square root of the matrix $S$.
Thus $A$ can be mapped to a positive definite symmetric matrix by the similarity transformation of the square root of $S$. That is, ${\tilde {A}}=\sqrt{S} A \sqrt{S}^{-1}$
 is positive definite symmetric.
Defined as $X=\sqrt{S}^{-1}Y$, we obtain the S-Oja-Brockett equation
$$\frac{dX}{dt}=AXB-XBX^TSAX.$$
Here $SA$ is symmetric.
This shows the S-Oja-Brockett equation is a real analytic gradient flow.
Thus we conclude {\rm (i)}. The proof of {\rm (ii)} is clear from the fact that 
the sorting properties of the S-Oja-Brockett equation clearly take over the properties of 
the Oja-Brockett equation.
{\rm (iii)} is obvious because the Oja-Brockett equation contains the Oja equation for the principal subspace as a special case.
{\hfill $\Box$}
}\end{pf}  

At the end of this section, we present simulation results for the S-Oja-Brockett equation with Eulerian discretization. 
\bigskip

Let $A=
\begin{pmatrix}
\frac{9}{2}     & -2              & 2   \\
 -2              & \frac{9}{2}    & 2   \\
  3               &    3            & 7      
\end{pmatrix}$,
$B=
\begin{pmatrix}
3     & 0      & 0   \\
0     & 2      & 0   \\
0     & 0      & 1      
\end{pmatrix}$ 
and
$S=
\begin{pmatrix}
  1               & 0                & 0   \\
  0               & 1                & 0   \\
  0               & 0                & \frac{2}{3}      
\end{pmatrix}.$ \\
Then $SA=
\begin{pmatrix}
\frac{9}{2}      & -2                & 2   \\
 -2               & \frac{9}{2}      & 2   \\
  2                & 2                 & \frac{14}{3}      
\end{pmatrix}.$
 \\  \bigskip

\begin{figure}[h]
    \centering
    \includegraphics[width=0.6\textwidth]{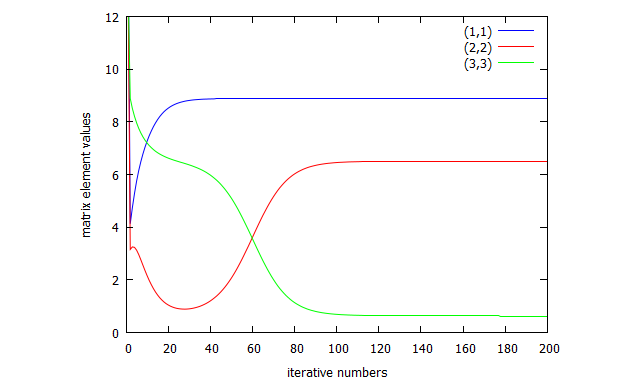}
    \caption{Diagonal components of $L(t)$}\label{figure_5}
\end{figure}

\begin{figure}[h]
    \centering
    \includegraphics[width=0.6\textwidth]{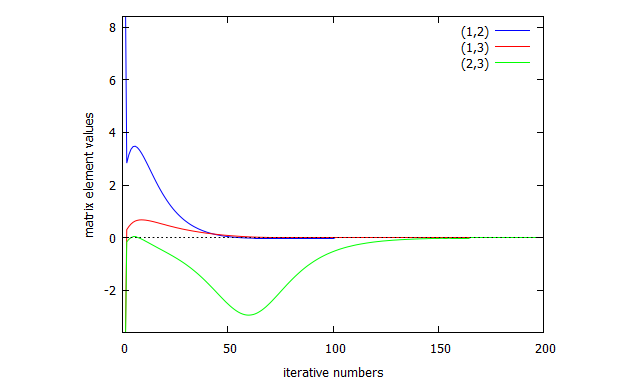}
    \caption{Off-diagonal components of $L(t)$}\label{figure_6} 
\end{figure}

The S-Oja-Brockett equation is discretized by the Euler method with a step size 0.01 such that 
$$
X_{n+1} = X_n + 0.01\cdot ( AX_n B - X_n BX_n^T SAX_n ), 
$$
and the following figures show the curves of the components of $L(t)=X(t)^TSAX(t)$
with an initial value  
$X_1=X(1)=
\begin{pmatrix}
  1               & 0                  & -1   \\
  0               & -1                & 1   \\
 -2               &-1                & 0      
\end{pmatrix}$.
The symbol $(i,j)$ in the figures represents the $i,j$ components of $L(t)$.
The eigenvalues of $A$ are obtained by LAPACK's routine, DGEEV as
$\{0.6193220895354239,\quad 6.500000000000001,\quad 8.880677910464575\},$
which means that the diagonal components of $L(t)$ converge to the eigenvalues of $A$. See Figure \ref{figure_5} and \ref{figure_6}.

  \qquad

\begin{figure}[h]
\begin{center}
\includegraphics[width=0.6\textwidth]{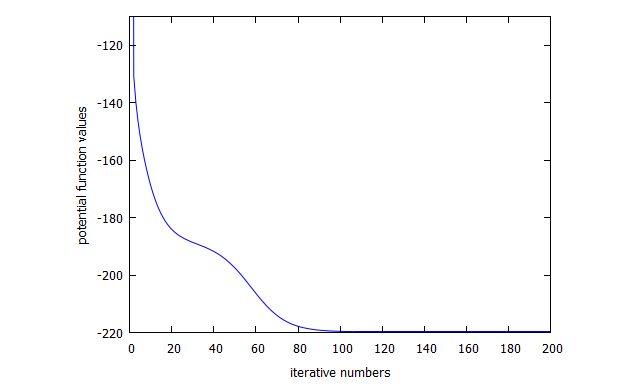}
\end{center}
\caption{Potential function values along the solution of the S-Oja-Brockett equation}
\label{figure_7}
\end{figure}

Figure \ref{figure_7} shows that the potential function descends monotonically 
along the solution of the S-Oja-Brockett equation.

\section{ S-Oja-Brockett equation for saddle point matrices}
In this section we investigate spectral properties of 
block matrices of the form
\begin{equation}\label{saddle} 
{\mathcal A}=
\begin{pmatrix}
P   & Q^T \\
-Q & R          
\end{pmatrix},
\end{equation}
where $P \in {\mathbb R}^{n \times n}$ is symmetric positive definite, 
$Q \in {\mathbb R}^{m \times n}$ has full rank with $m \leq n$, and
$R \in {\mathbb R}^{m \times m}$ is symmetric positive semidefinite. 
The matrix form (\ref{saddle}) can arise, for example, from finite element discretizations
of linearized Navier-Stokes equations and Maxwell equations, nonlinear optimization problems
and so on; see Benzi-Golub-Liesen \cite{BGL} and Benzi-Simoncini \cite{BS} for details.

In Section 3, the product of the axisymmetric elements of the matrix was non-negative, but in this section we consider the case where the product of the axisymmetric elements includes both non-negative and non-positive elements.
Specifically, consider the spectrum with respect to the following family of matrices:
\begin{align*}
{\mathfrak A}_{n,m}=
\biggl\{ \quad {\mathcal A} & = 
\begin{pmatrix}
P   & Q^T \\
-Q & R          
\end{pmatrix}
\in \mathbb{R}^{(n+m) \times (n+m)} \quad  \bigg| \quad 
0 < P = P^T \in \mathbb{R}^{n \times n},  \\ 
&  Q \in \mathbb{R}^{m \times n}, \quad m \leq n, \quad  {\rm rank}(Q)=m,  \quad  
0 \leq R = R^T \in \mathbb{R}^{m \times m}  \\
&  \quad {\rm s.t.} \quad  \biggl[   
\quad \exists  {\mathcal S} \in \mathbb{R}^{(n+m) \times (n+m)}
\quad  {\rm s.t.} \quad  0 < {\mathcal S} = {\mathcal S}^T, \quad  
{\mathcal A}^T{\mathcal S} ={\mathcal S}{\mathcal A} \quad
\biggr]
\quad \biggr\}.
\end{align*}

The objective is to concretize the conditions that ${\mathcal A}$ is diagonalizable and 
the eigenvalues of ${\mathcal A}$ are positive real numbers, 
and to find the eigenvalues and eigenvectors of ${\mathcal A}$ by the S-Oja-Brockett equation. 
Shen-Huang-Cheng \cite{SHC} proposed weaker conditions than those proposed by Liesen \cite{Lies}. 
However, the condition in \cite{SHC} is not a necessary and sufficient condition, 
thus we will show the necessary and sufficient conditions and give the necessary condition explicitly 
by the eigenvalues and singular values of the block matrix of ${\mathcal A}$.
Unlike the proof of \cite{SHC}, 
our proof does not require a case separation and slightly simplifies the relationship between $\varepsilon$ perturbations 
and the consisting of the eigenvalues and singular values of the block matrices.
\bigskip

We prepare some notations.
Let 
\begin{equation}\label{s_saddle} 
{\mathcal S}=
\begin{pmatrix}
P   & Q^T \\
Q & -R          
\end{pmatrix},
\end{equation}
where $P \in {\mathbb R}^{n \times n}$ is symmetric positive definite, 
$Q \in {\mathbb R}^{m \times n}$ has full rank with $m \leq n$, and
$R \in {\mathbb R}^{m \times m}$ is symmetric positive semidefinite. 
${\mathcal A}$ and ${\mathcal S}$ are perturbed by a scalar matrix, respectively such as 
\begin{equation}\label{d-saddle} 
{\mathcal A}_\delta =
\begin{pmatrix}
P +\delta I_n   & Q^T \\
-Q & \delta I_m + R          
\end{pmatrix},
\end{equation}
\begin{equation}\label{e-saddle} 
{\mathcal S}_\varepsilon =
\begin{pmatrix}
P -\varepsilon I_n   & Q^T \\
Q & \varepsilon I_m - R          
\end{pmatrix},
\end{equation}
where $0 \leq \delta$ and $0 < \varepsilon$.

Then we see the following  Sylvester equation.
\begin{lem}\label{perturb}
\begin{equation}\label{saddle_syl}
{\mathcal A}_\delta^T {\mathcal S}_\varepsilon =  {\mathcal S}_\varepsilon {\mathcal A}_\delta.
\end{equation}
\end{lem}  
\begin{pf}\rm{ By direct matrix computations.
{\hfill $\Box$}
}\end{pf}  

From Lemma \ref{perturb}, we have 
\begin{pro}
$\tilde{\mathcal A}_{\varepsilon, \delta}=
\begin{pmatrix}
{\mathcal A}_\delta  \\
{\mathcal S}_\varepsilon          
\end{pmatrix}
$ is a frame for a Lagrangian subspace of ${\mathbb R}^{2(n+m)}$.  
\end{pro}
Since we want to know about the spectral properties of a given matrix ${\mathcal A}$, 
we consider the frame, $\tilde{A}_\varepsilon=\tilde{A}_{\varepsilon,0}$ in the following discussion.
If there exists a positive definite ${\mathcal S}_\varepsilon$, 
then a given similarity transformation can make 
${\mathcal A}_\delta$ a symmetric matrix.
From the decomposition formulas for 
${\mathcal S}_\varepsilon$ with $\varepsilon \ne \lambda(P)$  and ${\mathcal A}={\mathcal A}_0$ 
such that 
\begin{align*}
{\mathcal S}_\varepsilon & =
\begin{pmatrix}
P -\varepsilon I_n   & Q^T \\
Q & \varepsilon I_m - R          
\end{pmatrix} \\
& =
\begin{pmatrix}
I_n                                &   O  \\
Q(P-\varepsilon I_n)^{-1}  &   I_m
\end{pmatrix}
\begin{pmatrix}
P -\varepsilon I_n  &   O  \\
O   & W_\varepsilon
\end{pmatrix}
\begin{pmatrix}
I_n   &    (P-\varepsilon I_n)^{-1}Q^T  \\
O    &     I_m
\end{pmatrix},
\end{align*}
where $W_\varepsilon = \varepsilon I_m - R - Q(P-\varepsilon I_n)^{-1}Q^T$
is the Schur complement of $P - \varepsilon I_n$ in ${\mathcal S}_\varepsilon$ and 
$$
{\mathcal A} =
\begin{pmatrix}
P    & Q^T \\
-Q  &  R          
\end{pmatrix}
=
\begin{pmatrix}
I_n          &   O  \\
-QP^{-1}  &   I_m
\end{pmatrix}
\begin{pmatrix}
P   &   O  \\
O   & W
\end{pmatrix}
\begin{pmatrix}
I_n   &     P^{-1}Q^T  \\
O    &     I_m
\end{pmatrix},
$$
where $W = R + QP^{-1}Q^T>0$, we obtain the following theorem:  
\bigskip

\begin{thm}[{\bf Positive definiteness of ${\mathcal S}_\varepsilon$}]\label{positive} 
Let ${\mathcal A} \in \mathfrak{A}_{n,m}$ be satisfying
$\lambda_{max}(R) < \lambda_{min}(P)$.
Then the necessary and sufficient condition for 
${\mathcal S}_\varepsilon$ to have positive eigenvalues are given by
{\rm (i),(ii)}, and {\rm (iii)}:
\begin{spacing}{1.5}
\leftline{$
{\rm (i)} \quad 0< P-\varepsilon I_n, \quad 
{\rm (ii)} \quad 0 < \varepsilon I_m -R, \quad 
{\rm (iii)} \quad Q(P-\varepsilon I_n)^{-1}Q^T < \varepsilon I_m - R.
$} 
\end{spacing}
Furthermore, 
using the singular value decomposition 
$Q=U\Sigma V^T$, where
$UU^T = I_m,  VV^T = I_n$ and  
$\Sigma=
\begin{pmatrix}
\Delta &  0 
\end{pmatrix}
 \in {\mathbb R}^{m \times n}$ with 
$\Delta = {\rm diag}(\sigma_1(Q),...,\sigma_m(Q))$,  
we have the sufficient condition for ${\mathcal S}_{\varepsilon}$ to be positive definite defined by 
{\rm (iv),(v)} and {\rm (vi)}:
\begin{spacing}{1.5}
\leftline{$
{\rm (iv)} \quad 0 < \lambda_{min}(P) - \varepsilon, 
\quad {\rm (v)} \quad  0 <  \varepsilon - \lambda_{max}(R), 
\quad  {\rm (vi)} \quad
\sigma_{max}^2(Q)<(\lambda_{min}(P) - \varepsilon)(\varepsilon - \lambda_{max}(R)). 
$} 
\end{spacing}

The condition under which $\varepsilon$ perturbations exist is given by 
\begin{spacing}{1.5}
\leftline{$
{\rm (vii)}\quad 2\sigma_{max}(Q) \leq \lambda_{min}(P) - \lambda_{max}(R).
$} 
\end{spacing}

Solving inequality {\rm (vi)} in $\varepsilon$ under condition {\rm (vii)}, we obtain
\begin{flushleft}
{\rm (viii)}$\quad \varepsilon_{-} < \varepsilon < \varepsilon_{+}, \quad$ where 
\end{flushleft}
\begin{align*}
\varepsilon_{\pm} & = 
\frac{\lambda_{min}(P) + \lambda_{max}(R)}{2} \\
& \qquad 
\pm \frac{1}{2}\sqrt{ ( \lambda_{min}(P) -2\sigma_{max}(Q) - \lambda_{max}(R) )
( \lambda_{min}(P) +2\sigma_{max}(Q) - \lambda_{max}(R) )} .
\end{align*}

\end{thm}
\begin{pf}\rm{
Conditions ${\rm (i), (ii),}$ and ${\rm (iii)}$ can be seen from the decomposition for 
${\mathcal S}_{\varepsilon}$. 
Substituting the singular value decomposition, 
$Q = U\Sigma V^T$
into ${\rm (iii)}$, we see 
\begin{equation}\label{positivein}
\Delta 
\begin{pmatrix}
I_m &  0 
\end{pmatrix}
(V^TPV -\varepsilon I_n)^{-1}
\begin{pmatrix}
I_m   \\
0
\end{pmatrix}
\Delta < \varepsilon I_m - U^TRU.
\end{equation}

From the inequality (\ref{positivein}), we have
\begin{equation}\label{lambda}
\lambda_{max} \left \{
\Delta 
\begin{pmatrix}
I_m & 0 
\end{pmatrix}
(V^TPV -\varepsilon I_n)^{-1}
\begin{pmatrix}
I_m   \\
0
\end{pmatrix}
\Delta \right \} < 
\lambda_{max}(\varepsilon I_m - U^TRU ).
\end{equation}
Since 
\begin{align*}
\lambda_{max} & \left \{
\begin{pmatrix}
\Delta & 0 
\end{pmatrix}
(V^TPV -\varepsilon I_n)^{-1}
\begin{pmatrix}
\Delta  \\
0
\end{pmatrix}
\right \} 
=
\lambda_{max} \left \{
\begin{pmatrix}
\Delta & 0 
\end{pmatrix}
\begin{pmatrix}
\Delta  \\
0
\end{pmatrix}
(V^TPV -\varepsilon I_n)^{-1}
\right \} 
\\ 
&= \lambda_{max} \left \{
\begin{pmatrix}
\Delta^2 & 0 \\
0 & 0  
\end{pmatrix}
(V^TPV -\varepsilon I_n)^{-1}
\right \}  = 
\frac{\sigma_{max}^2(Q)}{\lambda_{min}(V^TPV - \varepsilon)}     
= \frac{\sigma_{max}^2(Q)}{\lambda_{min}(P) - \varepsilon},
\end{align*}
substituting this into the inequality (\ref{lambda}),
we obtain the following inequality.
$$\frac{\sigma_{max}^2(Q)}{\lambda_{min}(P) - \varepsilon} < \varepsilon - \lambda_{max}(R).$$
This means {\rm (vi)}. 
{\rm (vii)} is a condition for the discriminant of the quadratic inequality 
{\rm (vi)} with respect to $\varepsilon$
 to have a real solution, and 
{\rm (viii)} is the solution of the quadratic inequality with respect to $\varepsilon$.
{\hfill $\Box$}
}\end{pf}  
\bigskip

Let ${\mathcal S}_{\varepsilon}$ have positive eigenvalues and 
let $\lambda_{\varepsilon}$ be an eigenvalue of ${\mathcal S}_{\varepsilon}$.
From Theorem \ref{positive} and Theorem 6 given in  
Dau\v{z}ickait\.{e}-Lawless-Scott-van Leeuwen \cite{DLSL}, 
we easily obtain the following fact.
\begin{pro}\label{eigen}
Suppose ${\mathcal S}_\varepsilon$ satisfies condition {\rm (iv),(v)} and {\rm(vi)} of Theorem \ref{positive}. 
Then the positive eigenvalues of ${\mathcal S}_{\varepsilon}$ lie in the interval 
$[\lambda_{-}(\varepsilon), \lambda_{+}(\varepsilon)]$,
where 
\begin{align*}
\lambda_{\pm}(\varepsilon) & =
\frac{\lambda_{min}(P) -\lambda_{max}(R)}{2} \\
& \qquad \pm 
\frac{1}{2}\sqrt{(\lambda_{max}(R) + \lambda_{min}(P))^2 + 4\sigma_{max}^2(Q) 
- 4\varepsilon \lambda_{min}(P) + \varepsilon (\lambda_{max}(R) -\varepsilon)  }
\end{align*}
\end{pro}
\bigskip

As an example, consider the following matrix given in Liesen \cite{Lies} and 
Shen-Huang-Cheng \cite{SHC} with $b=1/4$ and $c=1/12$.
\begin{eqnarray*}
{\mathcal A} = \left(\begin{array}{ccc|cc}
1               &  0     &  0  & 1/4    &  0       \\
0               &  2     &  0  &  0      &  1/4    \\
0               &  0     &  3  &  0      &  0       \\ \hline
-1/4 &  0     &  0  &  1/6   & -1/12  \\
0              & -1/4 &  0  & -1/12 &  1/6
\end{array}
\right).
\end{eqnarray*}
\bigskip

In this case, we see $\lambda_{min}(P)=1, \quad \sigma_{max}^2(Q)=1/16, \quad \lambda_{max}(R)=1/4$.
Thus, the sufficient conditions {\rm (iv)}, {\rm (v)}, and {\rm (vi)} are satisfied, and 
${\mathcal A}$ is positive. 
Since ${\mathcal A}$ is not a symmetric matrix, it is symmetrized by ${\mathcal S}_{\varepsilon}$, 
but it is not necessary to find $\varepsilon$ strictly. For such $\varepsilon$, 
since the minimum eigenvalue of $P$ is 1 and the maximum eigenvalue of $R$ is 1/4, 
we may set $\varepsilon=1/2$ and use the S-Oja-Brockett equation for 
${\mathcal A}$ to find eigenvalues. 
Note that the calculation of eigenvalues of ${\mathcal A}$ does not require 
the calculation of the square root of ${\mathcal S_{\varepsilon}}$, 
but the eigenvectors require the calculation of the square root of ${\mathcal S}_{\varepsilon}$.

To extract the three large eigenvalues and eigenvectors, 
let $X(t)$ be a $5 \times 3$ matrix.
To extract the five eigenvalues and eigenvectors, let $X(t)$ be a $5 \times 5$ matrix.
To find three and five eigencomponents, respectively, $B_3$and $B_5$ are defined as follows, 
and ${\mathcal S}_{1/2}$ is set as follows.
\bigskip

$B_3=
\begin{pmatrix}
3     &0        &0         \\
0      &2       &0         \\
0      &0        &1             
\end{pmatrix}.
\quad 
B_5=
\begin{pmatrix}
5     &0       &0       &0      &0      \\
0      &4      &0       &0       &0      \\
0      &0       &3      &0       &0      \\
0     &0        &0       &2       &0      \\
0     &0        &0       &0        &1   
\end{pmatrix}.
\quad 
{\mathcal S}_{1/2}=
\begin{pmatrix}
1/2     & 0      &0      &1/4       &0      \\
0     &3/2       &0      &0       &1/4      \\
0     &0       &5/2      &0       &0      \\
1/4     &0       &0      &1/3       &1/12      \\
0     &1/4       &0      &1/12       &1/3   
\end{pmatrix}.$
 \\  \bigskip

The five eigenvalues of ${\mathcal A}$ are obtain by a numerical routine, DGEEV in LAPACK such as 
$\{ 0.9153889054734862, 0.3188183561243709, 0.1339492200428448,  1.965176851692632,
 3.0 \}.$ The values of these eigenvalues are compared with the results of the following simulations.

\begin{figure}[h]
    \centering
    \includegraphics[keepaspectratio, scale=0.4]{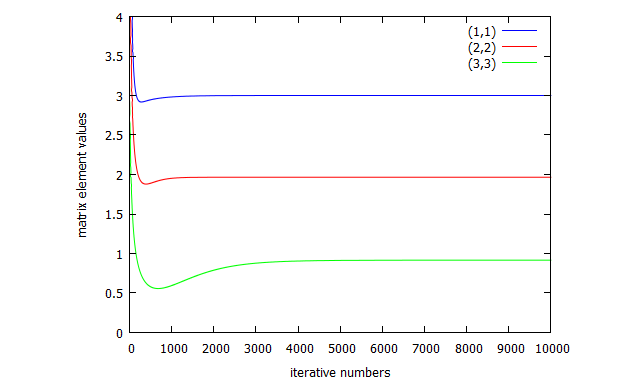}
    \caption{$L_{ii}$ components }\label{figure_8}
\end{figure}

\begin{figure}[h]
    \centering
    \includegraphics[keepaspectratio, scale=0.4]{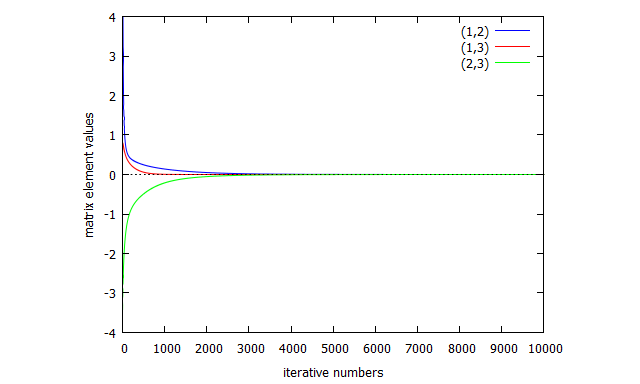}
    \caption{ Some $L_{ij}$ components $(i \ne j)$ }\label{figure_9} 
\end{figure}

Fig.\ref{figure_8} and Fig.\ref{figure_9} with $L(t)=X(t)^T{\mathcal S}_{1/2}{\mathcal A}X(t)$ show the results of simulating the S-Oja-Brockett equation with $B_3$ for eigenvalues of ${\mathcal A}$:
$$\frac{dX}{dt}={\mathcal A}XB_3-XB_3X^T{\mathcal S}_{1/2}{\mathcal A}X, \quad X \in {\mathbb R}^{5 \times 3}$$
and their corresponding eigenvectors with Eulerian differences:
$$
X_{n+1} = X_n + 0.001\cdot ( {\mathcal A}X_n B_3 - X_n B_3X_n^T {\mathcal S}_{1/2}{\mathcal A}X_n ),
\quad X_i \in {\mathbb R}^{5 \times 3}.
$$
Here, 
$X_1=
\begin{pmatrix}
1     & 0       &-1    \\
0     &-1      &1      \\
-2   &-1      &0      \\
0     &1       &1      \\
2     &-1      &0    
\end{pmatrix}.$
\bigskip


\begin{figure}[h]
\centering
\includegraphics[scale=0.4]{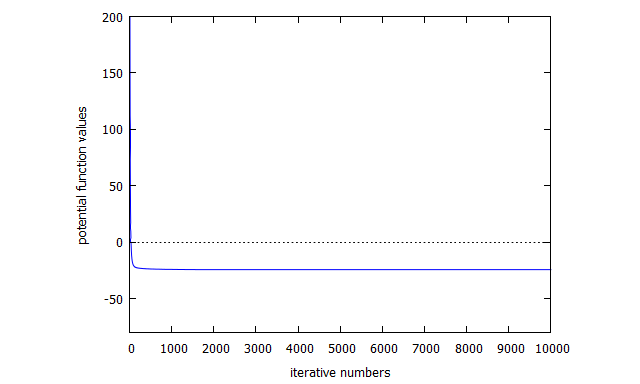}
\caption{Potential function values along the solution of the S-Oja-Brockett equation}
\label{figure_10}
\end{figure}

Figure \ref{figure_10} shows that the potential function descends monotonically 
along the solution of the S-Oja-Brockett equation with $B_3$.
From Fig.\ref{figure_8} and \ref{figure_9}, it can be seen that the three large eigenvalues obtained 
by LAPACK's routine, DGEEV, have been extracted.
Fig.\ref{figure_8} shows that the diagonal components of $L$ asymptotically approach the eigenvalues, 
and the certainty of this can be understood by the fact that the off-diagonal components of $L$ are approaching zero. See Figure \ref{figure_9}. From Fig. \ref{figure_10}, it can also be seen that the asymptotic approach to the optimal solution is achieved by monotonically descending the non-convex potential function.
\bigskip

\begin{figure}[h]
    \centering
    \includegraphics[keepaspectratio, scale=0.4]{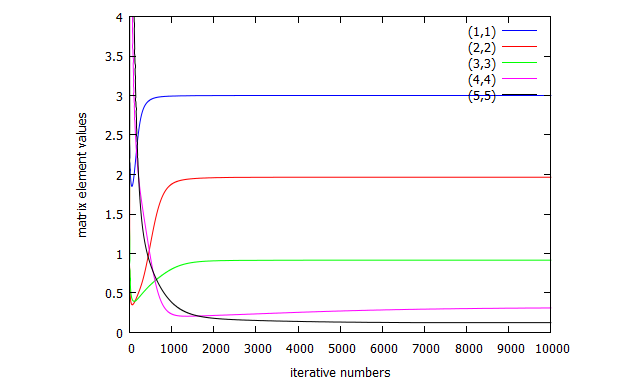}
    \caption{$L_{ii}$ components}\label{figure_11}
\end{figure}

\begin{figure}[h]
    \centering
    \includegraphics[keepaspectratio, scale=0.4]{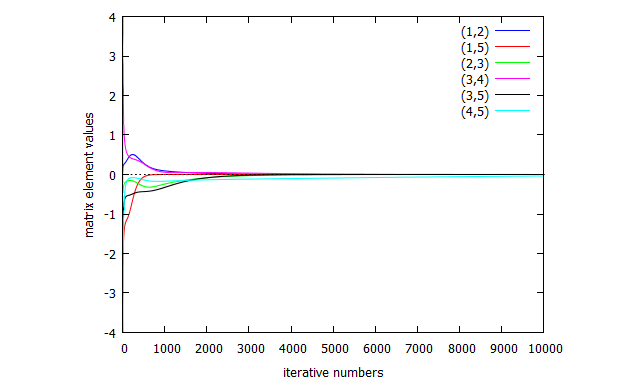}
    \caption{Some $L_{ij}$ components $(i \ne j)$}\label{figure_12} 
\end{figure}

Fig.\ref{figure_11} and Fig.\ref{figure_12} with $L(t)=X(t)^T{\mathcal S}_{1/2}{\mathcal A}X(t)$ show that the results of simulating the S-Oja-Brockett equation $B_5$ for eigenvalues of ${\mathcal A}$: 
$$\frac{dX}{dt}={\mathcal A}XB_5-XB_5X^T{\mathcal S}_{1/2}{\mathcal A}X, \quad X \in {\mathbb R}^{5 \times 5}$$
and their corresponding eigenvectors with Eulerian differences:
$$
X_{n+1} = X_n + 0.001\cdot ( {\mathcal A}X_n B_5 - X_n B_5X_n^T {\mathcal S}_{1/2}{\mathcal A}X_n ),
\quad X_i \in {\mathbb R}^{5 \times 5}.
$$
Here, 
$X_1=
\begin{pmatrix}
1     & 0      &-1      &2       &3      \\
0     &-1       &1      &6       &1      \\
-2     &-1       &0      &1       &9      \\
0     &-1       &1      &-3       &4      \\
-2     &-1       &0      &1       &3   
\end{pmatrix}.$
 \\  \bigskip
 

It can be seen from Fig.\ref{figure_11} that when all eigenvalues are obtained, 
a number of iterations are required, especially for convergence of small eigenvalues close to zero. See Figure \ref{figure_12}. 

\begin{figure}[h]
\centering
\includegraphics[scale=0.4]{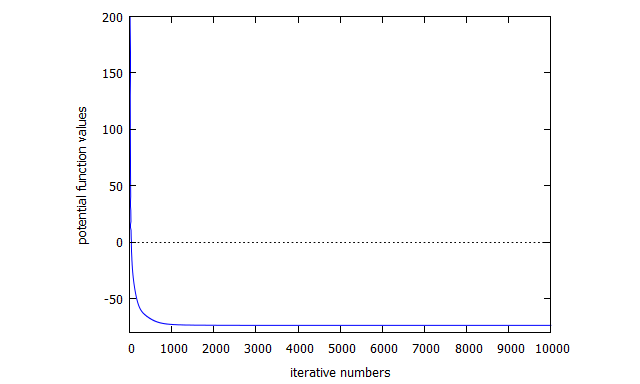}
\caption{Potential function values along the solution of the S-Oja-Brockett equation}\label{figure_13}
\end{figure}

Figure \ref{figure_13} shows that the potential function descends monotonically 
along the solution of the S-Oja-Brockett equation with $B_5$. 
Clearly lower descent than for the three main eigenvalues case.

\begin{figure}[h]
    \centering
    \includegraphics[keepaspectratio, scale=0.4]{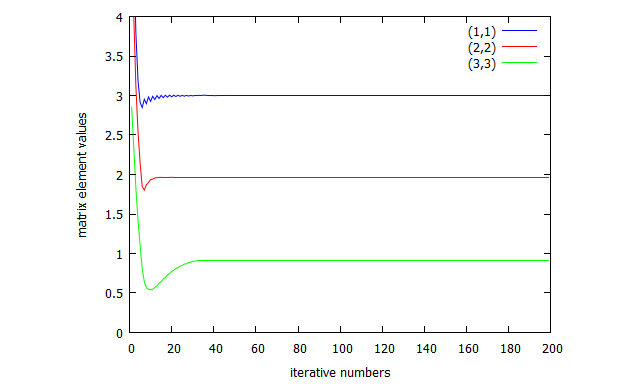}
    \caption{$L_{ii}$ components }\label{figure_14}
\end{figure}

\begin{figure}[h]
    \centering
    \includegraphics[keepaspectratio, scale=0.4]{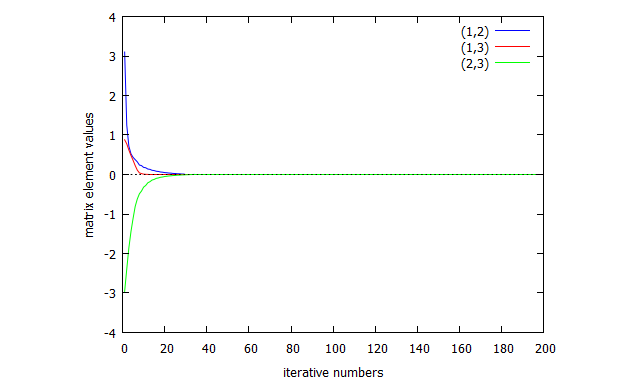}
    \caption{ Some $L_{ij}$ components $(i \ne j)$ }\label{figure_15} 
\end{figure}

Discretization of differential equations with a fixed step size requires changing the step size and checking convergence, depending on the size of the matrix and the presence of eigenvalues that are close together or close to zero. Therefore, we consider the following discretization that automatically determines the step size for each update of $X$.
\begin{thm}[{\bf Variable Eulerian discretization}]\label{VED} 
Let $A \in {\mathbb R}^{N \times N}$ be symmetrizable and have distinct positive eigenvalues.
And let $B \in {\mathbb R}^{M \times M}$ be a positive definite diagonal matrix with $M \leq N$.
For the S-Oja-Brockett equation
$$\frac{dX}{dt} = AXB-XBX^TSAX,$$
where $S (=S^T>0)$ is a symmetrizer of $A$, i.e., 
$\sqrt{S}A{\sqrt{S}}^{-1}$ is symmetric, 
we consider its Euler discretization defined by 
$$
X_{n+1} = X_n + \gamma(n) \cdot V(X_n), 
$$  
where $V(X_n) = AX_nB-X_nBX_n^TSAX_n$ and $0 < \gamma(n)$.
Then we obtain  
\begin{flushleft}
{\rm (i)} \quad The optimal step size $\gamma_{opt}(n)$ of the Euler discretization is given by 
\end{flushleft}
$$\gamma_{opt}(n) = {\rm min} \{ \quad \gamma \in {\mathbb R} \quad | \quad 0 < \gamma, \quad 
c_3(n) \cdot \gamma^3 + c_2(n) \cdot \gamma^2 + c_1(n) \cdot \gamma + c_0(n) = 0 \quad \},$$
\\
where 
\begin{align*}
c_3(n) & = {\rm tr}\{ (W_1(X_n))^2 \}, \\
c_2(n) & = 3{\rm tr}\{ W_1(X_n)W_2(X_n) \}, \\
c_1(n) & = 2{\rm tr}\{(W_2(X_n))^2 \} + {\rm tr}\{ W_1(X_n)W_3(X_n) \} + 2{\rm tr} \{W_4(X_n) \}, \\
c_0(n) & = {\rm tr}\{W_2(X_n)W_3(X_n))-{\rm tr} \{W_5(X_n)\}, \\
W_1(X_n) & = SAV(X_n)BV(X_n)^T, \\
W_2(X_n) & = SAX_nBV(X_n)^T, \\
W_3(X_n) & = SAX_nBX_n^T, \\
W_4(X_n) & = SA^2V(X_n)B^2V(X_n)^T, \\
W_5(X_n) & = SA^2V(X_n)B^2X_n^T.
\end{align*}
\begin{flushleft}
{\rm (ii)} \quad $\sqrt{S}X_{\infty}$ is the generalized eigenvectors of $A$. 
\end{flushleft}
\begin{flushleft}
{\rm (iii)} \quad $X_{\infty}^TSAX_{\infty}$ is a diagonal matrix. i.e., the generalized eigenvalues of $A$. 
\end{flushleft}
\end{thm}
\begin{pf} 
The S-Oja-Brockett equation is a gradient flow of the potential function 
$$g(X) = \frac{1}{4}{\rm tr}[(SAXBX^T)^2] - \frac{1}{2}{\rm tr}(SA^2XB^2X^T)$$
with respect to the Riemanian metric defined by 
$$\langle \Omega_1, \Omega_2 \rangle ={\rm tr}(SA \Omega_1 B \Omega_2^T), 
\quad \Omega_1,\Omega_2 \in {\mathbb R}^{N \times  M}.$$ 
That is, $\frac{dY}{dt}=-{\rm grad}g(X)$.
Therefore, with a one-step update from $X_n$ to $X_{n+1}$, 
we can find the step size optimization by finding the extreme value of 
$\phi(\gamma) = g(X_n + \gamma \cdot V(X_n))$. Since the cubic polynomial $\frac{d\phi}{d\gamma}$ with respect to $\gamma$,  
we obtain the optimal solution of  
$$\frac{d\phi}{d\gamma} = c_3(n) \cdot \gamma^3 + c_2(n) \cdot \gamma^2 + c_1(n) \cdot \gamma + c_0(n) = 0.$$ 
Let $\varphi(\gamma) = \frac{d\phi}{d\gamma}(\gamma)$. 
Since $c_3(n) >0$ for $X \ne 0$ and $\varphi (0) = c_0(n) = -{\rm tr}(SAV(X)BV(X)^T) < 0$ for $X \ne 0$, 
the cubic polynomial equation $\varphi(\gamma) = 0$ has at least one positive solution. 
Thus, we obtained {\rm (i)}.
\begin{flushleft}
{\rm (ii)} \quad By the Oja-Brockett equation, we know that
$$\sqrt{S}A\sqrt{S}^{-1} \cdot Y_{\infty} = Y_{\infty} \Lambda, $$
where $\Lambda$ is a diagonal matrix whose diagonal elements are the eigenvalues of $A$.
Substituting $Y_{\infty} = \sqrt{S}X_{\infty}$ into this, we see
$AX_{\infty} = X_{\infty} \Lambda.$
\end{flushleft}
\begin{flushleft}
{\rm (iii)} \quad By the Oja-Brockett equation, we see that
$Y_{\infty}^T\sqrt{S}A\sqrt{S}^{-1}Y_{\infty} = X_{\infty}^TSAX_{\infty}$ is a diagonal matrix.
\end{flushleft}
{\hfill $\Box$}
\end{pf}  
\bigskip
The results of the difference algorithm that automatically updates the Euler step size are again confirmed for  the example of Liesen \cite{Lies} and Shen-Huang-Cheng \cite{SHC}.
Fig.\ref{figure_14} and Fig.\ref{figure_15} with $L(t)=X(t)^T{\mathcal S}_{1/2}{\mathcal A}X(t)$ show the results of simulating the S-Oja-Brockett equation with $B_3$ for eigenvalues of ${\mathcal A}$:
$$\frac{dX}{dt}={\mathcal A}XB_3-XB_3X^T{\mathcal S}_{1/2}{\mathcal A}X, \quad X \in {\mathbb R}^{5 \times 3}$$
and their corresponding eigenvectors with Eulerian differences:
$$
X_{n+1} = X_n + \gamma(n)\cdot ( {\mathcal A}X_n B_3 - X_n B_3X_n^T {\mathcal S}_{1/2}{\mathcal A}X_n ),
\quad X_i \in {\mathbb R}^{5 \times 3}.
$$
Here, 
$X_1=
\begin{pmatrix}
1     & 0       &-1    \\
0     &-1      &1      \\
-2   &-1      &0      \\
0     &1       &1      \\
2     &-1      &0    
\end{pmatrix}.$


\begin{figure}[h]
\centering
\includegraphics[scale=0.4]{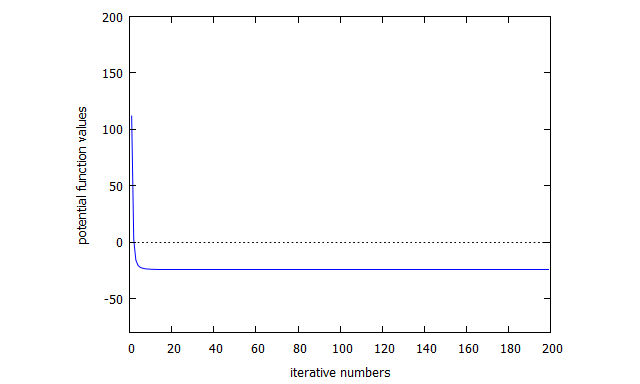}
\caption{Potential function values along the solution of the S-Oja-Brockett equation}\label{figure_16}
\end{figure}

Figure \ref{figure_16} shows that the potential function descends monotonically 
along the solution of the S-Oja-Brockett equation with $B_3$.

From Fig.\ref{figure_14} and \ref{figure_15}, it can be seen that the three large eigenvalues obtained 
by LAPACK's routine, DGEEV, have been extracted.
Fig.\ref{figure_14} shows that the diagonal components of $L$ asymptotically approach the eigenvalues, 
and the certainty of this can be understood by the fact that the off-diagonal components of $L$ are approaching zero. From Fig. \ref{figure_16}, it can also be seen that the asymptotic approach to the optimal solution is achieved by monotonically descending the non-convex potential function.
\bigskip

\begin{figure}[h]
    \centering
    \includegraphics[keepaspectratio, scale=0.4]{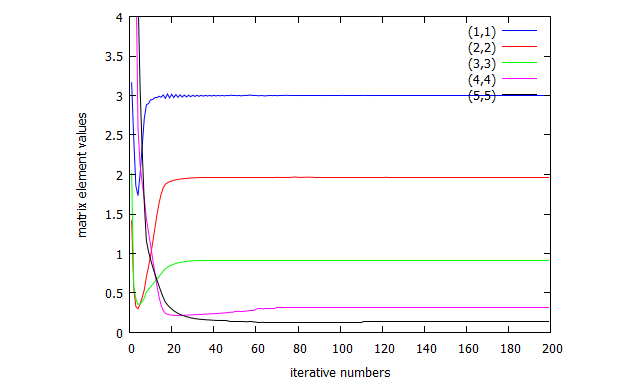}
    \caption{$L_{ii}$ components}\label{figure_17}
\end{figure}

\begin{figure}[h]
    \centering
    \includegraphics[keepaspectratio, scale=0.4]{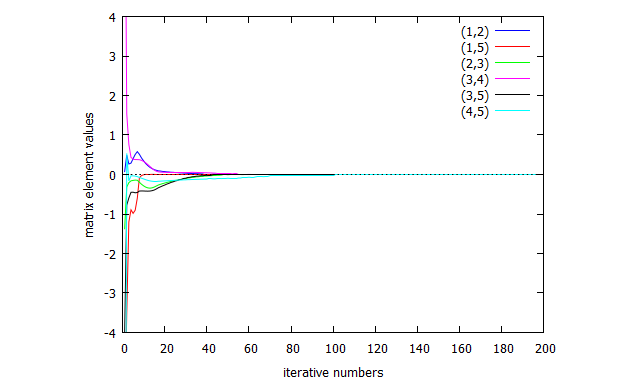}
    \caption{Some $L_{ij}$ components $(i \ne j)$}\label{figure_18} 
\end{figure}

Fig.\ref{figure_17} and Fig.\ref{figure_18} with $L(t)=X(t)^T{\mathcal S}_{1/2}{\mathcal A}X(t)$ show that the results of simulating the S-Oja-Brockett equations $B_5$ for eigenvalues of ${\mathcal A}$: 
$$\frac{dX}{dt}={\mathcal A}XB_5-XB_5X^T{\mathcal S}_{1/2}{\mathcal A}X, \quad X \in {\mathbb R}^{5 \times 5}$$
and their corresponding eigenvectors with Eulerian differences:
$$
X_{n+1} = X_n + \gamma(n)\cdot ( {\mathcal A}X_n B_5 - X_n B_5X_n^T {\mathcal S}_{1/2}{\mathcal A}X_n ),
\quad X_i \in {\mathbb R}^{5 \times 5}.
$$
Here, 
$X_1=
\begin{pmatrix}
1     & 0      &-1      &2       &3      \\
0     &-1       &1      &6       &1      \\
-2     &-1       &0      &1       &9      \\
0     &-1       &1      &-3       &4      \\
-2     &-1       &0      &1       &3   
\end{pmatrix}.$
 \\  \bigskip
 

\begin{figure}[h]
\centering
\includegraphics[scale=0.4]{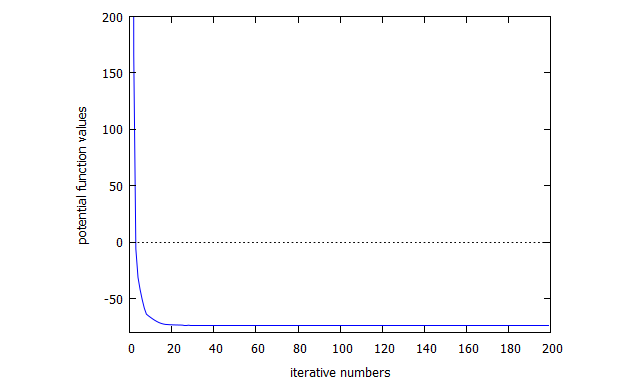}
\caption{Potential function values along the solution of the S-Oja-Brockett equation}\label{figure_19}
\end{figure}

It can be seen from Fig.\ref{figure_17} that when all eigenvalues are obtained, 
a number of iterations are required, especially for convergence of small eigenvalues close to zero.
See Figure \ref{figure_18}.
Figure \ref{figure_19} shows that the potential function descends monotonically 
along the solution of the S-Oja-Brockett equation with $B_5$. 
Clearly lower descent than for the three main eigenvalues.
Simulating variable Eulerian differences for an example by Liesen \cite{Lies} and 
Shen-Huang-Cheng \cite{SHC} 
with $b=1/4$ and $c=1/12$, the convergence is clearly about 50 times faster than Eulerian discretization with fixed step size.

The topics considered in this paper belong to a field where diverse disciplines such as neural networks, mathematical systems, numerical analysis, and real algebraic geometry intersect. The author hopes that the results obtained will contribute to the study of dynamical systems and numerical computation.

\section*{Acknowledgment}
This paper was inspired by Professor Yusaku Yamamoto's concrete comments 
on how to construct a diagonalized matrix $S$ for the matrix $A$ in Example 1 in Section 2.
The author would also like to thank Professor Uwe Helmke and Professor John B. Moore 
for the opportunity to discuss their book \cite{HM} and mathematical ideas with them in person more than 7 years ago.


\bigskip

\begin{thebibliography}{9}

\bibitem{Ama} S.-I. Amari.
Neural Theory of Association and Concept-Formation, 
Biol. Cybernetics 26 (1977) 175-185.

\bibitem{BGL} M. Benzi, G.H. Golub and J. Liesen, 
Numerical solution of saddle point problems, Acta Numerica 14 (2005) 1-137.

\bibitem{BS} M. Benzi and V. Simoncini, 
On the eigenvalues of a class of saddle point matrices, 
Numer. Math., 103 (2006) 173-196.

\bibitem{Bro1} R. W. Brockett, 
Dynamical systems that sort lists and solve linear programming problems, Proc. IEEE Conf. Decision and Control,
Austin, TX, (1988) 779-803. See also (Brockett, 1991).

\bibitem{Bro2} R. W. Brockett, 
Least squares matching problems, Linear Algebra
Appl. 122/123/124 (1989) 701-777.

\bibitem{Bro3} R. W. Brockett, 
Dynamical systems that sort lists, diagonalize matrices and solve linear programming problems, Linear Algebra Appl.
146 (1991) 79-91.

\bibitem{Chu} M. T. Chu, 
On the continuous realization of iterative processes,
SIAM Review 30 (1988) 375-387.

\bibitem{Cou} R. Courant, 
Zur Theorie der kleinen Schwingungen, Zeits. f. angew.
Math. Mech. 2 (1922) 278-285.

\bibitem{DLSL} I. Dau\v{z}ickait\.{e}, A. S. Lawless, J. A. Scott and P. J. van Leeuwen,
Spectral estimates for saddle point matrices arising in weak constraint four-dimensional variational data assimilation,
Numer. Linear Algebra Appl. Vol.27, Issue 5 (2020) e2313.

\bibitem{Fis} E. Fischer, 
Uber quadratische Formen mit reellen Koeffizienten, 
Monatsh. Math. u. Physik 16 (1905) 234-249.

\bibitem{HM} U.Helmke and J.B.Moore,
Optimization and dynamical systems, 
Springer-Verlag, London, 1996.

\bibitem{Hodge} J.H.Hodges,
Some matrix equations over finite field,
Ann. Mat. Pura Appl.(4), 44 (1957) 245-250.

\bibitem{Lies} J. Liesen, 
A note on the eigenvalues of saddle point matrices, Technical Report 10-2006, Institute of Mathematics, TU
 Berlin, 2006.

\bibitem{Lo} S. \L{}ojasiewicz,
Sur les trajectoires du gradient d’une fonction analytique. — Seminari
di Geometria, Bologna, Vol.15 (1983) 115-117. 

\bibitem{VonN} J. von Neumann, 
Some matrix-inequalities and metrization of matric-spaces, Tomsk Univ. Rev.1 (1937) 286-300.

\bibitem{Oja1} E. Oja,
A Simplified Neuron Model as a Principal Component Analyzer, 
J. Math. Biology 15 (1982) 267-273.

\bibitem{Oja2} E. Oja,
On Stochastic Approximation of the Eigenvectors and Eigenvalues of the Expectation of a Random Matrix,  
J. Mathematical Analysis and Applications 106 (1985) 69-84. 

\bibitem{Oja3} E. Oja,  
Neural networks, principal components and subspaces, 
Int. J. Neural Syst. 1(1) (1989) 61-68.

\bibitem{SHC} S-Q. Shen, T-Z. Huang and G-H. Cheng,
A condition for the nonsymmetric saddle point matrix being diagonalizable and having real and positive eigenvalues,
J.Comp. and Appl. Math. 220 (2008) 8-12. 

\bibitem{TS} K.Tanabe and M. Sagae, 
An Exact Cholesky Decomposition and the Generalized Inverse of the Variance-Covariance Matrix of the Multinomial Distribution, with Applications,
J. Royal Stat. Soc. Series B,Vol.54,1 (1992) 211-219.

\bibitem{TS2} K.Tanabe and M.Sagae,
Pivoting strategy for rank-one modification of LDM-like factorization,
Numer Algor 2 (1992) 137-153.

\bibitem{TZ} O. Taussky and H. Zassenhaus, 
On the similarity transformation between a matrix and its transpose, 
Pacific J. Math. 9 (1959) 893-896.

\bibitem{WE} J.L. Wyatt and I.M. Elfadel,
Time-domain solutions of Oja’s equations. 
Neural Comp., Vol.7, No.5 (1995) 915-922.

\bibitem{Xu} L. Xu, 
Least mean square error recognition principle for self organizing neural nets, 
Neural Netw., Vol.6, No.5 (1993) .627-648.

\bibitem{YHS} S.Yoshizawa, U.Helmke and K. Starkov,
Convergence analysis for principal component flows,
Int. J. Appl. Math. Comput. Sci., Vol.11, No.1 (2001) 223-236.
Corrections: ibid. Vol. 12, No.2 (2002) 299.



\end{thebibliography}
\end{document}